\documentclass{elsart3-1}

\usepackage{amssymb,stmaryrd,amsmath,amsfonts}

\usepackage[english,francais]{babel}

\newtheorem{theorem}{Theorem}[section]
\newtheorem{lemma}[theorem]{Lemma}
\newtheorem{e-proposition}[theorem]{Proposition}
\newtheorem{corollary}[theorem]{Corollary}
\newtheorem{e-definition}[theorem]{Definition\rm}

\newtheorem{theoreme}{Th\'eor\`eme}[section]

\newtheorem{proposition}[theoreme]{Proposition}

\renewcommand{\theequation}{\arabic{equation}}
\usepackage[all]{xy}

\newcommand{\N}{\mathbb{N}}

\newcommand{\Z}{\mathbb{Z}}

\def\og{\leavevmode\raise.3ex\hbox{$\scriptscriptstyle\langle\!\langle$~}}
\def\fg{\leavevmode\raise.3ex\hbox{~$\!\scriptscriptstyle\,\rangle\!\rangle$}}

\journal{the Acad\'emie des sciences}
\begin{document}
\centerline{}
\begin{frontmatter}


\selectlanguage{english}
\title{Motivic decompositions of projective homogeneous varieties and change of coefficients}


\selectlanguage{english}
\author[authorlabel1]{Charles De Clercq},
\ead{de.clercq.charles@gmail.com}

\address[authorlabel1]{Universit\'e Paris VI, Equipe Topologie et g\'eom\'etrie alg\'ebriques, 4 place Jussieu 75252 Paris, France}


\medskip
\begin{center}
{\small Received *****; accepted after revision +++++\\
Presented by £££££}
\end{center}

\begin{abstract}
\selectlanguage{english}
We prove that under some assumptions on an algebraic group $G$, indecomposable direct summands of the motive of a projective $G$-homogeneous variety with coefficients in $\mathbb{F}_p$ remain indecomposable if the ring of coefficients is any field of characteristic $p$. In particular for any projective $G$-homogeneous variety $X$, the decomposition of the motive of $X$ in a direct sum of indecomposable motives with coefficients in any finite field of characteristic $p$ corresponds to the decomposition of the motive of $X$ with coefficients in $\mathbb{F}_p$. We also construct a counterexample to this result in the case where $G$ is arbitrary.\\
{\it To cite this article: A.
Name1, A. Name2, C. R. Acad. Sci. Paris, Ser. I 340 (2005).}

\vskip 0.5\baselineskip

\selectlanguage{francais}

\noindent{\bf R\'esum\'e} \vskip 0.5\baselineskip \noindent

{\bf D\'ecompositions motiviques des vari\'et\'es projectives homog\`enes et changement des coefficients.}
Nous prouvons que sous certaines hypoth\`eses sur un groupe alg\'ebrique $G$, tout facteur direct ind\'ecomposable du motif associ\'e \`a une vari\'et\'e projective $G$-homog\`ene \`a coefficients dans $\mathbb{F}_p$ demeure ind\'ecomposable si l'anneau des coefficients est un corps de caract\'eristique $p$. En particulier pour toute vari\'et\'e projective $G$-homog\`ene $X$, la d\'ecomposition du motif de $X$ comme somme directe de motifs ind\'ecomposables \`a coefficients dans tout corps fini de caract\'eristique $p$ correspond \`a la d\'ecomposition du motif de $X$ \`a coefficients dans $\mathbb{F}_p$. Nous exhibons de plus un contre-exemple \`a ce r\'esultat dans le cas o\`u le groupe $G$ est quelconque.\\
{\it Pour citer cet article~: A. Name1, A. Name2, C. R. Acad. Sci.
Paris, Ser. I 340 (2005).}
\end{abstract}
\end{frontmatter}

\selectlanguage{english}
\section*{Introduction}

Let $F$ be a field, $\Lambda$ be a commutative ring, $CM(F;\Lambda)$ be the category of \emph{Grothendieck Chow motives} with coefficients in $\Lambda$, $G$ a semi-simple affine algebraic group and $X$ a projective $G$-homogeneous $F$-variety. The purpose of this note is to study the behaviour of the complete motivic decomposition (in a direct sum of indecomposable motives) of $X\in CM(F;\Lambda)$ when changing the ring of coefficients. In the first part we prove some very elementary results in non-commutative algebra and find sufficient conditions for the tensor product of two connected rings to be connected. In the second part we show that under some assumptions on $G$, indecomposable direct summands of $X$ in $CM(F;\mathbb{F}_p)$ remain indecomposable if the ring of coefficients is any field of characteristic $p$ (Theorem \ref{isoclass}), since these conditions hold for the \emph{reduced endomorphism ring} of indecomposable direct summands of $X$. In particular theorem \ref{isoclass} implies that the complete decomposition of the motive of $X$ with coefficients in any finite field of characteristic $p$ corresponds to the complete decomposition of the motive of $X$ with coefficients in $\mathbb{F}_p$. Finally we show that theorem \ref{isoclass} doesn't hold for arbitrary $G$ by producing a counterexample.

Let $\Lambda$ be a commutative ring. Given a field $F$, an $F$-variety will be understood as a separated scheme of finite type over $F$. Given such $\Lambda$ and an $F$-variety $X$, we can consider $CH_i(X;\Lambda)$, the Chow group of $i$-dimensional cycles on $X$ modulo rational equivalence with coefficients in $\Lambda$, defined as $CH_i(X)\otimes_{\Z} \Lambda$. These groups are the first step in the construction of the category $CM(F;\Lambda)$ of \emph{Grothendieck Chow motives} with coefficients in $\Lambda$. This category is constructed as the \emph{pseudo-abelian envelope} of the category $CR(F;\Lambda)$ of \emph{correspondences} with coefficients in $\Lambda$. Our main reference for the construction and the main properties of these categories is \cite{EKM}. For a field extension $E/F$ and any correspondence $\alpha\in CH(X\times Y;\Lambda)$ we denote by $\alpha_E$ the pull-back of $\alpha$ along the natural morphism $(X\times Y)_E\rightarrow X\times Y$. Considering a morphism of commutative rings $\varphi:\Lambda\longrightarrow \Lambda'$ we define the two following functors. The \emph{change of base field} functor is the additive functor $res_{E/F}:CM(F;\Lambda)\longrightarrow CM(E;\Lambda)$ which maps any summand $(X,\pi)[i]\in CM(F;\Lambda)$ to $(X_E,\pi_E)[i]$ and any morphism $\alpha:(X,\pi)[i]\rightarrow(Y,\rho)[j]$ to $\alpha_E$. The \emph{change of coefficents} functor is the additive functor $coeff_{\Lambda'/\Lambda}:CM(F;\Lambda)\longrightarrow CM(F;\Lambda')$ which maps any summand $(X,\pi)[i]$ to $(X,(id\otimes \varphi)(\pi))[i]$ and any morphism $\alpha:(X,\pi)[i]\rightarrow(Y,\rho)[j]$ to $(id\otimes \varphi)(\alpha)$.

\textbf{Acknowledgements  } I am very grateful to Nikita Karpenko for his suggestions and his support during this work. I also would like to thank Fran\c{c}ois Petit and Maksim Zhykhovich. Finally I am grateful to the referee for the useful remarks.

\section{On the tensor product of connected rings}

Recall that a ring $A$ is \emph{connected} if there are no idempotents in $A$ besides $0$ and $1$.

\begin{e-proposition}Let $A$ be a finite and connected ring. Then any element $a$ in $A$ is either nilpotent or invertible. The set $\mathcal{N}$ of nilpotent elements in $A$ is a two-sided and nilpotent ideal.
\end{e-proposition}

In order to prove Proposition 1.1 we will need the following elementary lemma.

\begin{lemma}\label{proj}Let $A$ be a finite ring. An appropriate power of any element $a$ of $A$ is idempotent.
\end{lemma}
Proof.~~For any $a\in A$, the set $\{a^n,~n\in \N\}$ is finite, hence there is a couple $(p,k)\in \N^2$ (with $k$ non-zero) such that $a^p=a^{p+k}$. The sequence $(a^n)_{n\geq p}$ is $k$-periodic and for example if $s$ is the lowest integer such that $p<sk$, $a^{sk}$ is idempotent.

Proof of Proposition 1.1.~~For any $a\in A$, an appropriate power of $a$ is an idempotent by lemma \ref{proj}. Since $A$ is connected, this power is either $0$ or $1$, that is to say $a$ is either nilpotent or invertible.

We now show that the set $\mathcal{N}$ of nilpotent elements in $A$ is a two-sided ideal. First if $a$ is nilpotent in $A$, then for any $b$ in $A$, $ab$ and $ba$ are not invertible, hence $ab$ and $ba$ belong to $\mathcal{N}$.

It remains to show that the sum of two nilpotent elements in $A$ is nilpotent. Setting $\nu$ for the number of nilpotent elements in $A$, we claim that for any sequence $a_1$,$\dots$, $a_{\nu}$ in $\mathcal{N}$, $a_1...a_{\nu}=0$. Indeed if $a_{\nu+1}$ is any nilpotent in $A$ the finite sequence $\Pi_1=a_1$, $\Pi_2=a_1a_2$,$\dots$, $\Pi_{\nu+1}=a_1a_2...a_{\nu+1}$ consists of nilpotents and by the pigeon-hole principle $\Pi_k=\Pi_s$, for some $k$ and $s$ satisfying $1\leq k<s\leq \nu+1$. Therefore $\Pi_s=\Pi_ka_{k+1}...a_s=\Pi_k$ which implies that $\Pi_k(1-a_{k+1}...a_s)=0$ and $\Pi_k=0$ since $1-a_{k+1}...a_s$ is invertible. With this in hand it is clear that for any $a$ and $b$ in $\mathcal{N}$, $(a+b)^{\nu}=0$. Furthermore $\mathcal{N}^{\nu}=0$ and $\mathcal{N}$ is nilpotent.

\begin{corollary}\label{connexion}Let $A$ be a finite and connected $\mathbb{F}_p$-algebra endowed with a ring morphism $\varphi:A\longrightarrow \mathbb{F}_p$. Then the set $\mathcal{N}$ of nilpotent elements in $A$ is precisely $ker(\varphi)$. Furthermore for any connected $\mathbb{F}_p$-algebra $E$, $A\otimes_{\mathbb{F}_p}E$ is connected.
\end{corollary}
Proof.~~For any $a\in \mathcal{N}$ and $n\in \N^{\ast}$ such that $a^n=0$, $0=\varphi(a^n)=\varphi(a)^n$, hence $a$ lies in the kernel of $\varphi$. On the other hand if $\varphi(a)=0$, $a$ is not invertible thus $a$ is nilpotent and $\mathcal{N}=ker(\varphi)$.
Since $\mathcal{N}$ is nilpotent, $\mathcal{N}\otimes E$ is also nilpotent. The sequence
\[\xymatrix{
0\ar[r]&\mathcal{N}\otimes E\ar[r]&A\otimes E\ar[r]^-{\psi}&E\ar[r]&0
 }\]
is exact and we want to show that any idempotent $P$ in $A\otimes_{\mathbb{F}_p}E$ is either $0$ or $1$. Since $E$ is connected, $\psi(P)$ is either $0$ or $1$. We may replace $P$ by $1-P$ and so assume that $P$ lies in the kernel of $\psi$, which implies that the idempotent $P$ is nilpotent, hence $P=0$.

\section{Application to motivic decompositions of projective homogeneous varieties}

For any semi-simple affine algebraic group $G$, the full subcategory of $CM(F;\Lambda)$ whose objects are finite direct sums of twists of direct summands of the motives of projective $G$-homogeneous $F$-varieties will be denoted $CM_G(F;\Lambda)$. We now use corollary \ref{connexion} to study how motivic decompositions in $CM_G(F;\Lambda)$ behave when extending the ring of coefficients. A pseudo-abelian category $\mathcal{C}$ satisfies the \emph{Krull-Schmidt principle} if the monoid $(\mathfrak{C},\oplus)$ is free, where $\mathfrak{C}$ denotes the set of the isomorphism classes of objects of $\mathcal{C}$.

In the sequel $\Lambda$ will be a connected ring and $X$ an $F$-variety. A field extension $E/F$ is a \emph{splitting field} of $X$ if the $E$-motive $X_E$ is isomorphic to a finite direct sum of twists of Tate motives. The $F$-variety $X$ is \emph{geometrically split} if $X$ splits over an extension of $F$, and $X$ satisfies the \emph{nilpotence principle}, if for any field extension $E/F$ the kernel of the morphism $res_{E/F}:End(M(X))\longrightarrow End(M(X_E))$ consists of nilpotents. Any projective homogeneous variety (under the action of a semi-simple affine algebraic group) is geometrically split and satisfies the nilpotence principle (see \cite{chermer}), therefore if $\Lambda$ is finite the Krull-Schmidt principle holds for $CM_G(F;\Lambda)$ by \cite[Corollary 3.6]{hhfba}, and we can serenely deal with motivic decompositions in $CM_G(F;\Lambda)$.

Let $G$ be a semi-simple affine adjoint algebraic group over $F$ and $p$ a prime. The absolute Galois group $Gal(F_{sep}/F)$ acts on the Dynkin diagram of $G$ and we say that $G$ is of \emph{inner type} if this action is trivial. By \cite{chermer} the subfield $F_G$ of $F_{sep}$ corresponding to the kernel of this action is a finite Galois extension of $F$, and we will say that $G$ is \emph{p-inner} if $[F_G:F]$ is a power of $p$. We now state the main result.
\begin{theorem}\label{isoclass}Let $G$ be a semi-simple affine adjoint $p$-inner algebraic group and $M\in CM_G(F;\mathbb{F}_p)$. Then for any field $L$ of characteristic $p$, $M$ is indecomposable if and only if $coeff_{L/\mathbb{F}_p}(M)$ is indecomposable.
\end{theorem}
If $X$ is geometrically split the image of any correspondence $\alpha\in CH_{\dim(X)}(X\times X;\Lambda)$ by the \emph{change of base field} functor $res_{E/F}$ to a splitting field $E/F$ of $X$ will be denoted $\overline{\alpha}$. The \emph{reduced endomorphism ring} of any direct summand $(X,\pi)$ is defined as $res_{E/F}(End_{CM(F;\Lambda)}((X,\pi)))$ and denoted by $\overline{End}((X,\pi))$.

Let $X$ be a complete and irreducible $F$-variety. The pull-back of the natural morphism $Spec(F(X))\times X\longrightarrow X\times X$ gives rise to $mult:CH_{dim(X)}(X\times X;\Lambda)\longrightarrow CH_0(X_{F(X)};\Lambda)\longrightarrow\Lambda$ (where the second map is the usual \emph{degree} morphism). For any correspondence $\alpha\in CH_{\dim(X)}(X\times X;\Lambda)$, $mult(\alpha)$ is called the \emph{multiplicity} of $\alpha$ and we say that a direct summand $(X,\pi)$ given by a projector $\pi\in CH_{\dim(X)}(X\times X;\Lambda)$ is \emph{upper} if $mult(\pi)=1$. If $(X,\pi)$ is an upper direct summand of a complete and irreducible $F$-variety, the multiplicity $mult:End_{CM(F;\Lambda)}((X,\pi))\longrightarrow \Lambda$ is a morphism of rings by \cite[Corollary 1.7]{anisot}.

\begin{proposition}\label{upper}Let $G$ be a semi-simple affine algebraic group and $M=(X,\pi)\in CM(F;\mathbb{F}_p)$ the upper direct summand of the motive of an irreducible and projective $G$-homogeneous $F$-variety. Then for any field $L$ of characteristic $p$, $M$ is indecomposable if and only if $coeff_{L/\mathbb{F}_p}(M)$ is indecomposable.
\end{proposition}
Proof.~~Since the change of coefficients functor is additive and maps any non-zero projector to a non-zero projector, it is clear that if $coeff_{L/\mathbb{F}_p}(M)$ is indecomposable, $M$ is also indecomposable.
Considering a splitting field $E$ of $X$, the reduced endomorphism ring $\overline{End}(M):=\overline{\pi}\circ\overline{End}(X)\circ \overline{\pi}$ is connected since $M$ is indecomposable and finite. Corollary \ref{connexion}, with $A=\overline{End}(M)$, $E=L$ and $\varphi=mult$ implies that $\overline{End}(M)\otimes L=\overline{End}(coeff_{L/\mathbb{F}_p}(M))$ is connected, therefore by the nilpotence principle $End(coeff_{L/\mathbb{F}_p}(M))$ is also connected, that is to say $coeff_{L/\mathbb{F}_p}(M)$ is indecomposable.

Proof of Theorem \ref{isoclass}.~~Recall that $G$ is a semi-simple affine adjoint $p$-inner algebraic group and consider a projective $G$-homogeneous $F$-variety $X$. By \cite[Theorem 1.1]{group}, any indecomposable direct summand $M$ of $X$ is a twist of the upper summand of the motive of an irreducible and projective $G$-homogeneous $F$-variety $Y$, thus we can apply proposition \ref{upper} to each indecomposable direct summand of $X$.

\textbf{Remark :} If $\Lambda$ is a finite, commutative and connected ring, complete motivic decompositions in $CM(F;\Lambda)$ remain complete when the coefficients are extended to the residue field of $\Lambda$ by \cite[Corollary 2.6]{vishikyagita}, hence the study of motivic decompositions in $CM_G(F;\Lambda)$, where $\Lambda$ is any finite connected ring whose residue field is of characteristic $p$, is reduced to the study motivic decompositions in $CM_G(F;\mathbb{F}_p)$.

We now produce a counterexample to Theorem \ref{isoclass} in the case where the algebraic group $G$ doesn't satisfy the needed assumptions. Let $L/F$ be a Galois extension of degree $3$. By \cite[Section 7]{chermer}, the endomorphism ring $End(M(Spec(L)))$ of the motive associated with the $F$-variety $Spec(L)$ with coefficients in $\mathbb{F}_2$ is the $\mathbb{F}_2$-algebra of $Gal(L/F)$, i.e. $\frac{\mathbb{F}_2[X]}{(X^3-1)}\simeq \mathbb{F}_2\times \mathbb{F}_4$, hence $M(Spec(L))=M\oplus N$, with $End(N)=\mathbb{F}_4$ and both $M$ and $N$ are indecomposable. Now $End(res_{\mathbb{F}_4/\mathbb{F}_2}(N))=\mathbb{F}_4\otimes\mathbb{F}_4$ is not connected since $1\otimes \alpha+\alpha\otimes 1$ is a non-trivial idempotent for any $\alpha\in \mathbb{F}_4\setminus\mathbb{F}_2$, hence $res_{\mathbb{F}_4/\mathbb{F}_2}(N)$ is decomposable.\\
Consider the $(PGL_2)_L$-homogeneous $L$-variety $\mathbb{P}^1_L$. The \emph{Weil restriction} $\mathcal{R}(\mathbb{P}^1_L)$ is a projective homogeneous $F$-variety under the action of the Weil restriction of $(PGL_2)_L$, and the minimal extension such that $\mathcal{R}((PGL_2)_L)$ is of inner type is $L$. By \cite[Example 4.8]{weil}, the motive with coefficients in $\mathbb{F}_2$ of $\mathcal{R}(\mathbb{P}^1_L)$ contains two twists of $Spec(L)$ as direct summands, therefore at least two indecomposable direct summands of $\mathcal{R}(\mathbb{P}^1_L)$ split off over $\mathbb{F}_4$.

\end{document}